\documentclass{article}[12pt]
\usepackage[french, english]{babel}
\usepackage{graphicx}
\usepackage{amsmath}
\usepackage{amsfonts}
\usepackage{amssymb}
\usepackage{amsthm}
\usepackage{hyperref}
\usepackage{mathtools}
\usepackage{comment}

\input xy
\xyoption{all}

\newtheorem{thm}{Th\'eor\`eme}[section]

\theoremstyle{definition}

\theoremstyle{remark}
\newtheorem*{rem}{Remarque}

\newcommand{\z}{\mathbb{Z}}

\newcommand{\im}{\mathrm{Im}}

\newcommand{\res}{\mathrm{Res}}

\newcommand{\ab}{\mathrm{ab}}
\newcommand{\br}{\mathrm{Br}\,}
\newcommand{\nr}{\mathrm{nr}}

\newcommand{\brnr}{{\mathrm{Br}_{\nr}}}
\newcommand{\al}{\mathrm{al}}

\newcommand{\brnral}{\mathrm{Br}_{\nr,\al}}

\newcommand{\spec}{\mathrm{Spec}\,}
\renewcommand{\cal}[1]{\mathcal{#1}}

\title{Le groupe de Brauer non ramifi\'e sur un corps global de caract\'eristique positive}
\author{Giancarlo Lucchini Arteche\\[5mm]
{\it\small D\'epartement de math\'ematique, universit\'e Paris-Sud}\\
{\it\small b\^atiment 425, 91405 Orsay cedex, France}\\
{\small giancarlo.lucchini@math.u-psud.fr}
}
\date{}

\begin{document}

\selectlanguage{french}
\maketitle

\begin{abstract}
En utilisant un th\'eor\`eme de Gabber sur les alt\'erations, on d\'emontre un r\'esultat d\'ecrivant la partie de torsion premi\`ere \`a $p$
du groupe de Brauer non ramifi\'e d'une vari\'et\'e $V$ lisse et g\'eom\'etriquement int\`egre sur un corps global de caract\'eristique
$p$ au moyen de l'\'evaluation des \'el\'ements de $\br V$ sur ses points locaux.\\

Mots cl\'es : Groupe de Brauer, corps global, alt\'eration.
\end{abstract}

\selectlanguage{english}
\begin{abstract}
{\bf The unramified Brauer group over a global field of positive characteristic.} Using a theorem of Gabber on alterations, we prove a result
describing the prime-to-$p$ torsion part of the unramified Brauer group of a smooth and geometrically integral variety $V$ over a global field
of characteristic $p$ by evaluating the elements of $\br V$ at its local points.\\

\end{abstract}

\selectlanguage{french}

\section{Introduction}
Dans \cite{HarariDuke}, Harari d\'emontre le r\'esultat suivant, lequel \'etait d\'ej\`a apparu (sous une forme quantitative plus pr\'ecise) dans le cas de
l'espace projectif chez Serre pour les \'el\'ements de $2$-torsion du groupe de Brauer (cf. \cite{SerreBrnr}) :

\begin{thm}[Th\'eor\`eme 2.1.1 de \cite{HarariDuke}]
Soient $K$ un corps de nombres et $X$ une $K$-vari\'et\'e g\'eom\'etriquement int\`egre, projective et lisse, dont on note $K(X)$ le corps de fonctions.
Soient $\alpha$ un \'el\'ement de $\br (K(X))$ qui n'est pas dans $\br X$ et $U$ un ouvert de Zariski non vide de $X$ tel que $\alpha\in\br U$. Alors, il
existe un ensemble infini (de densit\'e non nulle) $I$ de places $v$ de $K$ telles que, si l'on note $K_v$ le compl\'et\'e de $K$ par rapport \`a $v$,
la fl\`eche $U(K_v)\to\br K_v$ induite par $\alpha$ prenne une valeur non nulle.
\end{thm}

Il est bien connu qu'il existe un ensemble fini $S$ de places de $K$ tel que, si l'on note $\cal O_{K,S}$ l'anneau des $S$-entiers de $K$, il existe un
mod\`ele propre et lisse $\cal X$ de $X$ sur $\cal O_{K,S}$ avec un ouvert $\cal U$, dont la fibre g\'en\'erique correspond \`a $U$, tel que $\alpha$ se
rel\`eve en un \'el\'ement de $\br\cal U$. Le crit\`ere valuatif de
propr\'et\'e nous permet alors de relever tout $K_v$-point de $\cal X$ en un $\cal O_v$-point de $\cal X$ d\`es que $v\not\in S$. En particulier, on
voit que l'application $U(K_v)\to\br K_v$ induite par $\alpha$ se factorise par $\br\cal O_v=0$ d\`es que le $K_v$-point de $U$ se rel\`eve en un
$\cal O_v$-point de $\cal U$, ce qui nous dit que cette application atteint la valeur nulle pour presque toute place (i.e. pour toute place \`a l'exception
d'un nombre fini d'entre elles). De plus, on obtient une r\'eciproque du r\'esultat qui dit que pour tout \'el\'ement $\alpha$ de $\br X$ et pour presque
toute place $v$ de $K$ l'application $X(K_v)\to\br K_v$ induite par $\alpha$ est triviale. Si l'on consid\`ere alors une vari\'et\'e $V$ qui n'est pas propre,
ce r\'esultat appliqu\'e \`a une compactification lisse de $V$ nous permet de d\'ecrire le groupe de Brauer
non ramifi\'e $\brnr V$ d'apr\`es les th\'eor\`emes de puret\'e de Grothendieck (cf. \cite[Proposition 4.2.3]{ColliotSantaBarbara}).\\

Soit maintenant $K$ un corps global de caract\'eristique $p>0$, i.e. le corps de fonctions d'une courbe $C$ sur un corps fini. La preuve du th\'eor\`eme
de Harari est totalement adaptable \`a ce cas, sauf pour le fait que la preuve se base sur les th\'eor\`emes de puret\'e de Grothendieck, lesquels ne
sont valables a priori en caract\'eristique positive que sur la partie de torsion premi\`ere \`a $p$ du groupe de Brauer. On peut alors d\'emontrer
le r\'esultat suivant en suivant presque mot pour mot la preuve de Harari :

\begin{thm}\label{theoreme harari duke}
Soit $K$ un corps global de caract\'eristique $p>0$. Soit $X$ une $K$-vari\'et\'e propre, lisse et g\'eom\'etriquement int\`egre dont on note $K(X)$
le corps de fonctions. Soit $U$ un ouvert de $X$ et $\alpha\in\br U\subset\br(K(X))$ un \'el\'ement d'ordre $n$ premier \`a $p$ tel que
$\alpha\not\in\br X$. Alors il existe un ensemble infini (de densit\'e non nulle) $I$ de places $v$ de $K$ telles que
l'application $U(K_v)\to\br K_v$ induite par $\alpha$ soit non triviale.
\end{thm}

Si l'on consid\`ere maintenant une $K$-vari\'et\'e $V$ qui n'est pas propre, il n'est pas a priori \'evident que l'on puisse
utiliser ce r\'esultat pour d\'ecrire le groupe $\brnr V\{p'\}$, o\`u $\{p'\}$ veut dire que l'on prend les \'el\'ements de torsion premi\`ere \`a $p$,
car on ne dispose pas d'\'equivalent en caract\'eristique positive au th\'eor\`eme
d'Hironaka sur la r\'esolution des singularit\'es, donc l'existence d'une compactification lisse de $V$ n'est pas assur\'ee. Cependant, gr\^ace \`a un
r\'esultat de Gabber qui pr\'ecise le th\'eor\`eme de De Jong sur les alt\'erations (cf. \cite[Exp. X, Theorem 2.1]{ILO}), Borovoi, Demarche et
Harari ont r\'ecemment d\'emontr\'e le sens ``facile'' de cette description (cf. \cite[Proposition 4.2]{BDH}), i.e. que pour tout \'el\'ement $\alpha$
de $\brnr V\{p'\}$ et pour presque toute place $v$ de $K$ l'application $V(K_v)\to\br K_v$ induite par $\alpha$ est triviale.\\

En utilisant le m\^eme r\'esultat de Gabber, on peut montrer une version affaiblie de la r\'eciproque. Le but de ce texte est d'en donner la
d\'emonstration.

\section{R\'esultat principal}
Le r\'esultat est le suivant.

\begin{thm}\label{theoreme sans compactification lisse corps global}
Soient $K$ un corps global de caract\'eristique $p$ et $V$ une $K$-vari\'et\'e lisse et g\'eom\'etriquement int\`egre. Soit $\alpha\in\br V$ un
\'el\'ement d'ordre $n$ premier \`a $p$. Alors $\alpha$ appartient \`a $\brnr V$ si et seulement si, pour presque toute place $v$ de $K$,
l'application $V(L)\to\br L$ induite par $\alpha$ est nulle pour toute extension $L/K_v$ finie et purement ins\'eparable.
\end{thm}

\begin{rem}
L'auteur ignore si cette version affaiblie du r\'esultat est optimale ou pas, i.e. s'il suffit de regarder les $K_v$-points
et non pas les $L$-points pour toute extension $L/K_v$ purement ins\'eparable. La d\'emonstration actuelle ne laisse aucune
piste sur comment arriver \`a un tel r\'esultat.

Il est cependant important de remarquer que ce d\'efaut n'est en
g\'en\'eral pas g\^enant au moment des applications, comme on le verra \`a la fin du texte.
\end{rem}

\begin{proof}
L'un des sens ayant d\'ej\`a \'et\'e d\'emontr\'e, il suffit de se concentrer sur le cas o\`u $\alpha\not\in\brnr V$. Il faut alors trouver une
infinit\'e de places $v$ de $K$, des extensions purement ins\'eparables $L/K_v$ et des $L$-points $P$ de $V$ tels que $\alpha(P)\neq 0$.
On note aussi que, puisque $\br V$ est un groupe de torsion, on peut supposer que $\alpha$ est d'ordre $\ell^m$ pour un certain $\ell\neq p$
premier et $m\geq 1$. 

Le th\'eor\`eme de Nagata (cf. \cite{ConradNagata}) nous dit qu'il existe une $K$-compactification $X$ de $V$ qui n'est pas forc\'ement lisse, i.e.
une $K$-vari\'et\'e propre $X$ munie d'une $K$-immersion ouverte $V\hookrightarrow X$.
Le th\'eor\`eme de Gabber nous dit alors qu'il existe une extension finie de corps $K'/K$ de degr\'e premier \`a $\ell$ et une $\ell'$-alt\'eration
$X'\to X_{K'}$, i.e. un morphisme propre, surjectif et g\'en\'eriquement fini de degr\'e premier \`a $\ell$, avec $X'$ une $K'$-vari\'et\'e lisse. On en
d\'eduit que $X'$ est propre et lisse sur $K'$. De plus, on peut supposer que $X'$ est g\'eom\'etriquement int\`egre. En effet, quitte \`a prendre une
composante connexe de $X'$ dominant $X$, on peut supposer qu'elle est connexe. Il suffit alors de noter que la factorisation de Stein assure l'existence
d'une factorisation $X'\to \spec K'' \to \spec K'$ telle que $X'$ est g\'eom\'etriquement connexe sur $K''$, donc quitte \`a changer $K'$ par $K''$
(extension qui reste de degr\'e premier \`a $\ell$ car $[K'':K']$ divise $[K'(X'):K'(X)]$), on a que $X'$ est g\'eom\'etriquement int\`egre sur $K'$ car elle
est lisse et g\'eome\'triquement connexe. On a alors le diagramme commutatif
suivant, o\`u les carr\'es sont cart\'esiens et les fl\`eches $\hookrightarrow$ repr\'esentent des immersions ouvertes :
\[\xymatrix{
V' \ar@{^{(}->}[r] \ar[d] & X' \ar[d]  \\
V_{K'} \ar@{^{(}->}[r] \ar[d] & X_{K'} \ar[d] \\
V \ar@{^{(}->}[r] & X 
}
\]
Puisque $\alpha\not\in\brnr V$, on sait qu'il existe un anneau de valuation discr\`ete $A$ de corps de fractions $K(V)$ et corps r\'esiduel
$\kappa\supset K$ tel que l'image de $\alpha$ par l'application r\'esidu
\[H^2(K(V),\mu_{\ell^m})\to H^1(\kappa,\z/\ell^m\z),\]
est non nulle, cf. par exemple \cite[\S2]{ColliotSantaBarbara}. Soit $K_0/K(V)$ la sous-extension s\'eparable maximale de $K'(V')$. En
consid\'erant la fermeture int\'egrale de $A$ dans $K_0$ et en localisant en un id\'eal premier convenable on trouve un anneau de valuation
discr\`ete $A_0\subset K_0$ contenant $A$, de corps de fractions $K_0$, et tel que le produit
$e_{A_0/A}f_{A_0/A}$ soit premier \`a $\ell$, o\`u $e$ et $f$ repr\'esentent l'indice de ramification et le degr\'e r\'esiduel respectivement.
Si l'on d\'efinit alors $A'$ comme la fermeture int\'egrale de $A_0$ dans $K'(V')$, on a que $A'$ est un anneau de valuation discr\`ete et que
$e_{A'/A_0}f_{A'/A_0}$ est une puissance de $p$, cf. la preuve de \cite[Ch. 4, Proposition 1.31]{Liu}.  Ainsi, on voit que $A'/A$ est une extension
d'anneaux de valuation discr\`ete avec $e_{A'/A}f_{A'/A}$ premier \`a $\ell$. De plus, on voit bien que $K'(V')$ est le corps de fractions de $A'$ et que,
si l'on note $\kappa'$ le corps r\'esiduel de $A'$, on a $\kappa'\supset K'$. Alors, d'apr\`es \cite[Proposition 3.3.1]{ColliotSantaBarbara}, on a le
diagramme commutatif
\[\xymatrix{
H^2(K(V),\mu_{\ell^m}) \ar[r] \ar[d]^{\res} & H^1(\kappa,\z/\ell^m\z) \ar[d]^{e_{A'/A}\cdot\res}\\
H^2(K'(V'),\mu_{\ell^m}) \ar[r] & H^1(\kappa',\z/\ell^m\z). \\
}\]
Soit $\kappa_1$ la sous-extension s\'eparable maximale de $\kappa'/\kappa$. On a alors la factorisation suivante de la fl\`eche verticale de droite :
\[H^1(\kappa,\z/\ell^m\z)\xrightarrow{e_{A'/A}\cdot\res} H^1(\kappa_1,\z/\ell^m\z)\xrightarrow{\res} H^1(\kappa',\z/\ell^m\z).\]
Puisque l'extension $\kappa'/\kappa_1$ est purement ins\'eparable, si l'on note $\Gamma_{\kappa_1}$ et $\Gamma_{\kappa'}$ les groupes de Galois
absolus correspondants, il est facile de v\'erifier que le morphisme canonique $\Gamma_{\kappa'}\to\Gamma_{\kappa_1}$ est un isomorphisme, ce
qui nous dit que la fl\`eche de droite est aussi un isomorphisme, tandis que pour la
fl\`eche de gauche on a que $e_{A'/A}[\kappa_1:\kappa]$ divise $e_{A'/A}f_{A'/A}$ et alors le morphisme est injectif par l'argument classique de
restriction-corestriction. En particulier, on voit que l'image $\alpha'\in \br V'$ de $\alpha$ n'est pas dans $\brnr V'\{p'\}=\br X'\{p'\}$ car l'image
de $\alpha'$ dans $H^1(\kappa',\z/\ell^m\z)$ est non nulle.\\

Le th\'eor\`eme \ref{theoreme harari duke} nous donne enfin un ensemble infini $I'$ de places $w$ de $K'$ telles qu'il existe un $K'_w$-point $P'_w$
de $V'$ tel que $\alpha'(P'_w)\neq 0$. Ce $K'_w$-point induit \'evidemment un $K'_w$-point $Q'_w$ de $V_{K'}$, donc de $V$, tel que
$\alpha(Q'_w)\neq 0$. Consid\'erons alors la sous-extension s\'eparable maximale $K\subset K_1\subset K'$. L'ensemble $I'$, qui est de densit\'e non
nulle, induit par restriction un ensemble $I_1$ de places de $K_1$ de densit\'e non nulle. D'autre part, on sait que la densit\'e de l'ensemble
$I_2$ des places de $K_1$ ayant m\^eme corps r\'esiduel que la place qu'elles induisent sur $K$ est \'egale \`a 1, cf. par exemple
\cite[p. 215]{Heilbronn}. On voit alors que le sous-ensemble $I$ de $\Omega_K$ form\'e des places induites par celles dans
$I_1\cap I_2\subset\Omega_{K_1}$ est de densit\'e non nulle, donc infini. De plus, quitte \`a enlever les places ramifi\'ees pour l'extension $K_1/K$
(qui forment un ensemble fini), on peut supposer que l'extension $K'_w/K_v$ est purement ins\'eparable car elle est alors d\'ecomposable en une suite
d'extensions radicielles de degr\'e $p$, tout comme $K'/K_1$. Les $K'_w$-points $Q'_w$ conviennent alors.
\end{proof}

\section{Applications}
Avec le th\'eor\`eme \ref{theoreme sans compactification lisse corps global}, on peut obtenir des versions en caract\'eristique positive de beaucoup
d'autres r\'esultats autour du groupe de Brauer non ramifi\'e des vari\'et\'es sur un corps de nombres. En voici un exemple pour les espaces
homog\`enes qui est l'analogue d'un r\'esultat de Harari (cf. \cite[Proposition 4]{HarariBulletinSMF}).

\begin{thm}\label{theoreme harari SMF}
Soient $K$ un corps global de caract\'eristique $p>0$ et $G$ un $K$-groupe alg\'ebrique fini d'ordre premier \`a $p$ plong\'e dans un $K$-groupe $G'$
semi-simple simplement connexe (par exemple
$G'=\mathrm{SL}_n$) et soit $V=G'/G$. Soient $G^\ab$ l'ab\'elianis\'e de $G$ et $M=\hat G^\ab$ son groupe des caract\`eres. Alors le groupe de
Brauer non ramifi\'e alg\'ebrique $\brnral V$ s'identifie au sous-groupe de $H^1(K,M)$ constitu\'e des \'el\'ements $\alpha$ ayant la propri\'et\'e
suivante :

Pour presque toute
place $v$ de $K$, la restriction $\alpha_v$ de $\alpha$ dans $H^1(K_v,M)$ est orthogonale (pour l'accouplement local donn\'e par le cup-produit)
au sous-ensemble $\im [H^1(K_v,G)\to H^1(K_v,G^\ab)]$ de $H^1(K,G^\ab)$.
\end{thm}

En outre, on peut avec le th\'eor\`eme \ref{theoreme sans compactification lisse corps global} enlever les hypoth\`eses d'existence d'une
compactification lisse dans certains r\'esultats de Borovoi, Demarche et Harari, comme par exemple leur th\'eor\`eme 7.4.b dans \cite{BDH}, en
rempla\c cant bien entendu le groupe de Brauer de la compactification par le groupe de Brauer non ramifi\'e.\\

Les d\'etails de ces r\'esultats, notamment la preuve des th\'eor\`emes \ref{theoreme harari duke} et \ref{theoreme harari SMF}, appara\^\i tront
dans la th\`ese de doctorat de l'auteur.

\section*{Remerciements}
L'auteur tient \`a remercier David Harari pour son support constant pendant la r\'edaction de cet article.

\bibliographystyle{plain}
\bibliography{Brnr}

\begin{thebibliography}{1}

\bibitem{BDH}
Mikhail Borovoi, Cyril Demarche, and David Harari.
\newblock Complexes de groupes de type multiplicatif et groupe de {B}rauer non
  ramifi\'e des espaces homog\`enes.
\newblock {P}reprint : http://arxiv.org/abs/1203.5964, {\`a} para{\^\i}tre dans
  a {A}nn. {S}ci. {E}c. {N}orm. {S}up., 2012.

\bibitem{ColliotSantaBarbara}
J.-L. Colliot-Th{\'e}l{\`e}ne.
\newblock Birational invariants, purity and the {G}ersten conjecture.
\newblock In {\em {$K$}-theory and algebraic geometry: connections with
  quadratic forms and division algebras ({S}anta {B}arbara, {CA}, 1992)},
  volume~58 of {\em Proc. Sympos. Pure Math.}, pages 1--64. Amer. Math. Soc.,
  Providence, RI, 1995.

\bibitem{ConradNagata}
Brian Conrad.
\newblock Deligne's notes on {N}agata compactifications.
\newblock {\em J. Ramanujan Math. Soc.}, 22(3):205--257, 2007.

\bibitem{HarariDuke}
David Harari.
\newblock M\'ethode des fibrations et obstruction de {M}anin.
\newblock {\em Duke Math. J.}, 75(1):221--260, 1994.

\bibitem{HarariBulletinSMF}
David Harari.
\newblock Quelques propri\'et\'es d'approximation reli\'ees \`a la cohomologie
  galoisienne d'un groupe alg\'ebrique fini.
\newblock {\em Bull. Soc. Math. France}, 135(4):549--564, 2007.

\bibitem{Heilbronn}
H.~Heilbronn.
\newblock Zeta-functions and {$L$}-functions.
\newblock In {\em Algebraic {N}umber {T}heory ({P}roc. {I}nstructional {C}onf.,
  {B}righton, 1965)}, pages 204--230. Thompson, Washington, D.C., 1967.

\bibitem{ILO}
Luc Illusie, Yves Laszlo, and Fabrice Orgogozo.
\newblock Travaux de {G}abber sur l'uniformisation locale et la cohomologie
  {\'e}tale des sch{\'e}mas quasi-excellents. {S}{\'e}minaire a l'{\'e}cole
  polytechnique 2006--2008.
\newblock Preprint : http://arxiv.org/abs/1207.3648, 2012.

\bibitem{Liu}
Qing Liu.
\newblock {\em Algebraic geometry and arithmetic curves}, volume~6 of {\em
  Oxford Graduate Texts in Mathematics}.
\newblock Oxford University Press, Oxford, 2002.
\newblock Translated from the French by Reinie Ern{\'e}, Oxford Science
  Publications.

\bibitem{SerreBrnr}
Jean-Pierre Serre.
\newblock Sp\'ecialisation des \'el\'ements de {${\rm Br}_2({\bf
  Q}(T_1,\cdots,T_n))$}.
\newblock {\em C. R. Acad. Sci. Paris S\'er. I Math.}, 311(7):397--402, 1990.

\end{thebibliography}

\end{document}